\documentclass{ws-ijnt}

\begin{document}

\markboth{Aritram Dhar}{On ${}_5\psi_5$ Identities of Bailey}

%
\catchline{}{}{}{}{}
%

\title{ON ${}_5\psi_5$ IDENTITIES OF BAILEY}

\author{Aritram Dhar}

\address{Department of Mathematics\\ University of Florida\\ 358 Little Hall,
Gainesville, FL 32611, USA\\
\email{aritramdhar@ufl.edu} }

\maketitle

\begin{history}
\received{22 March 2023}
\accepted{26 August 2023}
\end{history}

\begin{abstract}
In this paper, we provide proofs of two ${}_5\psi_5$ summation formulas of Bailey using a ${}_5\phi_4$ identity of Carlitz. We show that in the limiting case, the two ${}_5\psi_5$ identities give rise to two ${}_3\psi_3$ summation formulas of Bailey. Finally, we prove the two ${}_3\psi_3$ identities using a technique initially used by Ismail to prove Ramanujan's ${}_1\psi_1$ summation formula and later by Ismail and Askey to prove Bailey's very-well-poised ${}_6\psi_6$ sum.
\end{abstract}



\section{Introduction}\label{s1}	

Let $a$ and $q$ be variables and define the conventional $q$-Pochammer symbol $$(a)_n=(a;q)_n:=\prod_{k=0}^{n-1}(1-aq^k)$$ for any positive integer $n$ and $(a)_0=1$. For $\lvert q\rvert<1$, we define $$(a)_{\infty}=(a;q)_{\infty}:=\lim_{n\rightarrow\infty}(a;q)_n.$$ We define $(a)_n$ for all real numbers $n$ by $$(a)_n := \dfrac{(a)_{\infty}}{(aq^n)_{\infty}}.$$ For variables $a_1,a_2,\ldots,a_k$, we define the shorthand notations $$(a_1,a_2,\ldots,a_k;q)_n:=\prod_{i=1}^{k}(a_i;q)_n\, ,$$ $$(a_1,a_2,\ldots,a_k;q)_{\infty}:=\prod_{i=1}^{k}(a_i;q)_{\infty}.$$\par Next, we require the following formulas from Gasper and Rahman \cite[Appendix I]{5}
\begin{equation}\label{eq11}
(a;q)_{n+k}=(a;q)_n(aq^n;q)_k,    
\end{equation}
\begin{equation}\label{eq12}
(a;q)_{-n}=\dfrac{1}{(aq^{-n};q)_n}=\dfrac{(-q/a)^n}{(q/a;q)_n}q^{\binom{n}{2}},
\end{equation} 
\begin{equation}\label{eq13}
(aq^{-n};q)_k=\dfrac{(a;q)_k(q/a;q)_n}{(q^{1-k}/a;q)_n}q^{-nk},\quad \text{and}    
\end{equation}
\begin{equation}\label{eq14}
\dfrac{(a;q)_{n-k}}{(b;q)_{n-k}}=\dfrac{(a;q)_n}{(b;q)_n}\dfrac{(q^{1-n}/b;q)_k}{(q^{1-n}/a;q)_k}\left(\dfrac{b}{a}\right)^k.    
\end{equation}
\\\par We invite the reader to examine Gasper and Rahman’s text \cite{5} for an introduction to basic hypergeometric series, whose notations we follow. For instance, the ${}_r\phi_{r-1}$ unilateral and ${}_r\psi_r$ bilateral basic hypergeometric series with base $q$ and argument $z$ are defined, respectively, by
\begin{align*}
{}_r\phi_{r-1} \left[
	\setlength\arraycolsep{2pt}
	\begin{matrix}
		a_1,\ldots,a_r \\
		\multicolumn{2}{c}{
			\begin{matrix}
				b_1,\ldots,b_{r-1} 	
			\end{matrix}}
	\end{matrix} \hspace{1pt}
;q, z \right]&:=\sum_{k=0}^{\infty}\dfrac{(a_1,\ldots,a_r;q)_k}{(q,b_1,\ldots,b_{r-1};q)_k}z^k,\quad \lvert z\rvert<1, \\
\quad{}_r\psi_r \left[
   \begin{matrix}
      a_1,\ldots,a_r \\
      b_1,\ldots,b_r
   \end{matrix}
;q, z \right]&:=\sum_{k=-\infty}^{\infty}\dfrac{(a_1,\ldots,a_r;q)_k}{(b_1,\ldots,b_r;q)_k}z^k\,\,\,\,,\quad \left\lvert\dfrac{b_1\ldots b_r}{a_1\ldots a_r}\right\rvert<\lvert z\rvert<1. 
\end{align*}
\par Throughout the remainder of this paper, we assume that $\lvert q\rvert<1$. We now present the statements of the main identities which we prove in this paper.
\\
\begin{theorem}(Bailey \cite[eq. $3.1$]{2})\label{5psi51}
For any non-negative integer $n$,
\begin{equation}\label{eq15}
{}_5\psi_5 \left[
   \begin{matrix}
      b, &c, &d, &e, &q^{-n} \\
      q/b, &q/c, &q/d, &q/e, &q^{n+1}
   \end{matrix}
;q, q \right]=\dfrac{(q,q/bc,q/bd,q/cd;q)_n}{(q/b,q/c,q/d,q/bcd;q)_n}    
\end{equation}
where $bcde=q^{n+1}$.
\end{theorem}

\begin{theorem}(Bailey \cite[eq. $3.2$]{2})\label{5psi52}
For any non-negative integer $n$,
\begin{equation}\label{eq16}
{}_5\psi_5 \left[
   \begin{matrix}
      b, &c, &d, &e, &q^{-n} \\
      q^2/b, &q^2/c, &q^2/d, &q^2/e, &q^{n+2}
   \end{matrix}
;q, q \right]=\dfrac{(1-q)(q^2,q^2/bc,q^2/bd,q^2/cd;q)_n}{(q^2/b,q^2/c,q^2/d,q^2/bcd;q)_n}
\end{equation}
where $bcde=q^{n+3}$.
\end{theorem}

\begin{theorem}(Bailey \cite[eq. $2.2$]{2})\label{3psi31}
\begin{equation}\label{eq17}
{}_3\psi_3 \left[
   \begin{matrix}
      b, &c, &d \\
      q/b, &q/c, &q/d
   \end{matrix}
;q, \dfrac{q}{bcd} \right] = \dfrac{(q,q/bc,q/bd,q/cd;q)_{\infty}}{(q/b,q/c,q/d,q/bcd;q)_{\infty}}.
\end{equation}
\end{theorem}

\begin{theorem}(Bailey \cite[eq. $2.3$]{2})\label{3psi32}
\begin{equation}\label{eq18}
{}_3\psi_3 \left[
   \begin{matrix}
      b, &c, &d \\
      q^2/b, &q^2/c, &q^2/d
   \end{matrix}
;q, \dfrac{q^2}{bcd} \right] = \dfrac{(q,q^2/bc,q^2/bd,q^2/cd;q)_{\infty}}{(q^2/b,q^2/c,q^2/d,q^2/bcd;q)_{\infty}}.
\end{equation}
\end{theorem}
\bigskip 
Bailey \cite{2} proved Theorems \ref{3psi31} and \ref{3psi32} by letting $a\rightarrow 1$ and setting $a=q$ in the ${}_6\phi_5$ summation formula \cite[II.$20$]{5} respectively and mentioned that (\ref{eq15}) and (\ref{eq16}) follow from Jackson's $q$-analogue of Dougall's theorem \cite[II.$22$]{5}.\par Our work is motivated by Ismail's initial proof \cite{6} of Ramanujan's ${}_1\psi_1$ summation formula which can be stated as
\begin{equation}\label{eq19}
{}_1\psi_1 \left[
   \begin{matrix}
      a \\
      b
   \end{matrix}
;q, z \right]=\dfrac{(q,b/a,az,q/az;q)_{\infty}}{(b,q/a,z,b/az;q)_{\infty}}
\end{equation} where $\lvert b/a\rvert<\lvert z\rvert<1$ and by Askey and Ismail's proof \cite{1} of Bailey's very-well-poised ${}_6\psi_6$ identity which is
\begin{align*}
{}_6\psi_6 \left[
   \begin{matrix}
      q\sqrt{a}, &-q\sqrt{a}, &b, &c, &d, &e \\
      \sqrt{a}, &-\sqrt{a}, &aq/b, &aq/c, &aq/d, &aq/e
   \end{matrix}
;q, \dfrac{qa^2}{bcde} \right]\\ = \dfrac{(aq,aq/bc,aq/bd,aq/be,aq/cd,aq/ce,aq/de,q,q/a;q)_{\infty}}{(aq/b,aq/c,aq/d,aq/e,q/b,q/c,q/d,q/e,qa^2/bcde;q)_{\infty}}\numberthis\label{eq110}
\end{align*}
provided $\lvert qa^2/bcde\rvert<1$.\par To prove (\ref{eq19}) and (\ref{eq110}), Ismail \cite{6} and Askey and Ismail \cite{1} show that the two sides of (\ref{eq19}) and (\ref{eq110}) are analytic functions that agree infinitely often near a point that is an interior point of the domain of analyticity and hence they are identically equal.\\\par
To this end, we employ the following $q$-hypergeometric series identities
\\
\begin{theorem}(Carlitz \cite[eq. $3.4$]{3})\label{5phi4}
For any non-negative integer $n$,
\begin{align*}
{}_5\phi_4 \left[
	\setlength\arraycolsep{6pt}
	\begin{matrix}
		q^{-n}, &\quad\quad b, &\quad\quad c, &\quad\quad d, &e \\
		\multicolumn{5}{c}{
			\begin{matrix}
				q^{-n+1}/b, &q^{-n+1}/c, &q^{-n+1}/d, &q^{-n+1}/e 	
			\end{matrix}}
	\end{matrix} \hspace{2pt}
;q, q \right] \\ = q^{m(1+m-n)}(de)^{-m}\dfrac{(q^{-n})_{2m}(q^{-n+1}/bc,q^{-n+1}/bd,q^{-n+1}/be;q)_m}{(q,q^{-n+1}/b,q^{-n+1}/d,q^{-n+1}/e,q^{n-m}c;q)_m}(q^{2m-n})_{n-2m}\numberthis\label{eq111}
\end{align*}
where $m=\lfloor n/2\rfloor$ \textit{and} $bcde=q^{1+m-2n}$.
\end{theorem}
We note that for $n$ even, Theorem \ref{5phi4} is Chu's \cite[p.~$279$]{4} Corollary $3$ where $\delta=0$ and for $n$ odd, Theorem \ref{5phi4} is Chu's \cite[p.~$280$]{4} Corollary $7$ where $\delta=0$.
\\

\begin{theorem}(Jackson's terminating q-analogue of Dixon's sum \cite[II.$15$]{5})\label{3phi2Jackson}
For any non-negative integer $m$,
\begin{equation}\label{eq112}
{}_3\phi_2 \left[
	\setlength\arraycolsep{6pt}
	\begin{matrix}
		q^{-2m}, &\quad a, &\quad b \\
		\multicolumn{3}{c}{
			\begin{matrix}
				q^{-2m+1}/a, &q^{-2m+1}/b 	
			\end{matrix}}
	\end{matrix} \hspace{2pt}
;q, \dfrac{q^{-m+2}}{ab} \right] \\ = \dfrac{(a,b;q)_m(q,ab;q)_{2m}}{(q,ab;q)_m(a,b;q)_{2m}}.
\end{equation}
\end{theorem}

\begin{theorem}(Carlitz \cite[eq. $2.5$]{3})\label{3phi2Carlitz}
For any non-negative integer $n$,
\begin{align*}
{}_3\phi_2 \left[
	\setlength\arraycolsep{6pt}
	\begin{matrix}
		q^{-n}, &\quad a, &\quad b \\
		\multicolumn{3}{c}{
			\begin{matrix}
				q^{-n+1}/a, &q^{-n+1}/b 	
			\end{matrix}}
	\end{matrix} \hspace{2pt}
;q, \dfrac{q^{-n+m+1}z}{ab} \right] \\ = \sum_{2j\le n}(-1)^j\dfrac{(q^{-n})_{2j}(q^{-n+1}/ab)_j}{(q,q^{-n+1}/a,q^{-n+1}/b;q)_j}q^{-j(j-1)/2+mj}z^j(z)_{m-j}(q^{j+m-n}z)_{n-m-j}\numberthis\label{eq113}
\end{align*}
where $m=\lfloor n/2\rfloor$.    
\end{theorem}
\bigskip
The paper is organized as follows. In Section \ref{s2}, we give the proofs of the two ${}_5\psi_5$ identities (\ref{eq15}) and (\ref{eq16}) respectively. In Section \ref{s3}, we show that the two ${}_5\psi_5$ identities (\ref{eq15}) and (\ref{eq16}) become the two ${}_3\psi_3$ identities (\ref{eq17}) and (\ref{eq18}) respectively when $n\rightarrow \infty$. Finally, we provide proofs of the two ${}_3\psi_3$ identities (\ref{eq17}) and (\ref{eq18}) in Section \ref{s4}.

\section{Proofs of the Two ${}_5\psi_5$ Identities}\label{s2}

\subsection{Proof of Theorem \ref{5psi51}}
\begin{proof}
Replacing $n$ by $2m$, $b$ by $bq^{-m}$, $c$ by $cq^{-m}$, $d$ by $dq^{-m}$ and $e$ by $eq^{-m}$ in (\ref{eq111}), we get
\begin{align*}
{}_5\phi_4 \left[
	\setlength\arraycolsep{5pt}
	\begin{matrix}
		q^{-2m}, &\quad bq^{-m}, &\quad cq^{-m}, &\quad dq^{-m}, &\quad eq^{-m} \\
		\multicolumn{5}{c}{
			\begin{matrix}
				q^{-m+1}/b, &q^{-m+1}/c, &q^{-m+1}/d, &q^{-m+1}/e   
			\end{matrix}}
	\end{matrix} \hspace{2pt}
;q,q \right] \\ = q^{m^2+m}(de)^{-m}\dfrac{(q^{-2m})_{2m}(q/bc,q/bd,q/be;q)_m}{(q,q^{-m+1}/b,q^{-m+1}/d,q^{-m+1}/e,c;q)_m}\numberthis\label{eq21}
\end{align*} where $bcde=q^{m+1}$. Now, we have
\[{}_5\psi_5 \left[
   \begin{matrix}
      b, &c, &d, &e, &q^{-n} \\
      q/b, &q/c, &q/d, &q/e, &q^{n+1}
   \end{matrix}
;q, q \right]\]
\begin{align*}&=\sum_{k=-\infty}^{\infty}\dfrac{(b,c,d,e,q^{-n};q)_k}{(q/b,q/c,q/d,q/e,q^{n+1};q)_k}q^k\\
&=\sum_{k=-n}^{\infty}\dfrac{(b,c,d,e,q^{-n};q)_k}{(q/b,q/c,q/d,q/e,q^{n+1};q)_k}q^k\quad (\text{since}\,1/(q^{n+1})_k=0\,\text{for all}\,k<-n)\\
&=\sum_{k=0}^{\infty}\dfrac{(b,c,d,e,q^{-n};q)_{k-n}}{(q/b,q/c,q/d,q/e,q^{n+1};q)_{k-n}}q^{k-n}\\
&=\dfrac{(b,c,d,e,q^{-n};q)_{-n}q^{-n}}{(q/b,q/c,q/d,q/e,q^{n+1};q)_{-n}}\sum_{k=0}^{\infty}\dfrac{(q^{-2n},bq^{-n},cq^{-n},dq^{-n},eq^{-n};q)_k}{(q,q^{-n+1}/b,q^{-n+1}/c,q^{-n+1}/d,q^{-n+1}/e;q)_k}q^k\\
&=\dfrac{(b,c,d,e,q^{-n};q)_{-n}(q^{-2n})_{2n}(q/bc,q/bd,q/be;q)_nq^{n^2}}{(q/b,q/c,q/d,q/e,q^{n+1};q)_{-n}(q,q^{-n+1}/b,q^{-n+1}/d,q^{-n+1}/e,c;q)_n(de)^n}
\end{align*}
\\
where the last equality above follows from (\ref{eq21}) (after replacing $m$ by $n$). Then simplifying the last expression above using (\ref{eq11}), (\ref{eq12}) and (\ref{eq13}) with appropriate substitutions, we get
\\
\[{}_5\psi_5 \left[
   \begin{matrix}
      b, &c, &d, &e, &q^{-n} \\
      q/b, &q/c, &q/d, &q/e, &q^{n+1}
   \end{matrix}
;q, q \right] = \dfrac{(q,q/bc,q/bd,q/cd;q)_n}{(q/b,q/c,q/d,q/bcd;q)_n}\]
\\
where $bcde=q^{n+1}$ for $n\in \mathbb{N}\cup\{0\}$. This completes the proof of Theorem \ref{5psi51}.
\end{proof}
\bigskip

\subsection{Proof of Theorem \ref{5psi52}}
\begin{proof}
Replacing $n$ by $2m+1$, $b$ by $bq^{-m-1}$, $c$ by $cq^{-m-1}$, $d$ by $dq^{-m-1}$ and $e$ by $eq^{-m-1}$ in (\ref{eq111}), we get
\begin{align*}
{}_5\phi_4 \left[
	\setlength\arraycolsep{5pt}
	\begin{matrix}
		q^{-2m-1}, &bq^{-m-1}, &cq^{-m-1}, &dq^{-m-1}, &eq^{-m-1} \\
		\multicolumn{5}{c}{
			\begin{matrix}
				q^{-m+1}/b, &q^{-m+1}/c, &q^{-m+1}/d, &q^{-m+1}/e   
			\end{matrix}}
	\end{matrix} \hspace{2pt}
;q,q \right] \\ = (q-1)q^{m^2+2m-1}(de)^{-m}\dfrac{(q^{-2m-1})_{2m}(q^2/bc,q^2/bd,q^2/be;q)_m}{(q,q^{-m+1}/b,q^{-m+1}/d,q^{-m+1}/e,c;q)_m}.\numberthis\label{eq22}
\end{align*} where $bcde=q^{m+3}$. Now, we have
\[{}_5\psi_5 \left[
   \begin{matrix}
      b, &c, &d, &e, &q^{-n} \\
      q^2/b, &q^2/c, &q^2/d, &q^2/e, &q^{n+2}
   \end{matrix}
;q, q \right]\]
\begin{align*}&=\sum_{k=-\infty}^{\infty}\dfrac{(b,c,d,e,q^{-n};q)_k}{(q^2/b,q^2/c,q^2/d,q^2/e,q^{n+2};q)_k}q^k\\
&=\sum_{k=-n-1}^{\infty}\dfrac{(b,c,d,e,q^{-n};q)_k}{(q^2/b,q^2/c,q^2/d,q^2/e,q^{n+2};q)_k}q^k\quad (\text{since}\,1/(q^{n+2})_k=0\,\text{for all}\,k<-n-1)\\
&=\sum_{k=0}^{\infty}\dfrac{(b,c,d,e,q^{-n};q)_{k-n-1}}{(q^2/b,q^2/c,q^2/d,q^2/e,q^{n+2};q)_{k-n-1}}q^{k-n-1}\\
&=\dfrac{(b,c,d,e,q^{-n};q)_{-n-1}q^{-n-1}}{(q^2/b,q^2/c,q^2/d,q^2/e,q^{n+2};q)_{-n-1}}\sum_{k=0}^{\infty}\dfrac{(q^{-2n-1},bq^{-n-1},cq^{-n-1},dq^{-n-1},eq^{-n-1};q)_k}{(q,q^{-n+1}/b,q^{-n+1}/c,q^{-n+1}/d,q^{-n+1}/e;q)_k}q^k\\
&=\dfrac{(q-1)(b,c,d,e,q^{-n};q)_{-n-1}(q^{-2n-1})_{2n}(q^2/bc,q^2/bd,q^2/be;q)_nq^{n^2+n-2}}{(q^2/b,q^2/c,q^2/d,q^2/e,q^{n+2};q)_{-n-1}(q,q^{-n+1}/b,q^{-n+1}/d,q^{-n+1}/e,c;q)_n(de)^n}
\end{align*}
\\
where the last equality above follows from (\ref{eq22}) (after replacing $m$ by $n$). Then simplifying the last expression above using (\ref{eq11}), (\ref{eq12}) and (\ref{eq13}) with appropriate substitutions, we get
\\
\[{}_5\psi_5 \left[
   \begin{matrix}
      b, &c, &d, &e, &q^{-n} \\
      q^2/b, &q^2/c, &q^2/d, &q^2/e, &q^{n+2}
   \end{matrix}
;q, q \right]=\dfrac{(1-q)(q^2,q^2/bc,q^2/bd,q^2/cd;q)_n}{(q^2/b,q^2/c,q^2/d,q^2/bcd;q)_n}\]
\\
where $bcde=q^{n+3}$ for $n\in \mathbb{N}\cup\{0\}$. This completes the proof of Theorem \ref{5psi52}.    
\end{proof}

\section{Two Limiting Cases}\label{s3}
Letting $n\rightarrow \infty$ in (\ref{eq15}) and simplifying using (\ref{eq13}) with appropriate substitutions, we get
\[{}_3\psi_3 \left[
   \begin{matrix}
      b, &c, &d \\
      q/b, &q/c, &q/d
   \end{matrix}
;q, \dfrac{q}{bcd} \right]=\dfrac{(q,q/bc,q/bd,q/cd;q)_\infty}{(q/b,q/c,q/d,q/bcd;q)_\infty}\]
which is exactly (\ref{eq17}).
\\\par
Similarly, letting $n\rightarrow \infty$ in (\ref{eq16}) and simplifying using (\ref{eq13}) with appropriate substitutions, we get 
\[{}_3\psi_3 \left[
   \begin{matrix}
      b, &c, &d \\
      q^2/b, &q^2/c, &q^2/d
   \end{matrix}
;q, \dfrac{q^2}{bcd} \right]=\dfrac{(q,q^2/bc,q^2/bd,q^2/cd;q)_\infty}{(q^2/b,q^2/c,q^2/d,q^2/bcd;q)_\infty}\]
which is exactly (\ref{eq18}).

\section{Ismail Type Proofs of the Two ${}_3\psi_3$ Identities}\label{s4}
In this Section, we derive the the two ${}_3\psi_3$ identities (\ref{eq17}) and (\ref{eq18}) using Ismail's method \cite{6}.
\subsection{Proof of Theorem \ref{3psi31}}
\begin{proof}
Replacing $a$ by $bq^{-m}$ and $b$ by $cq^{-m}$ in (\ref{eq112}), we get
\begin{equation}\label{eq41}
{}_3\phi_2 \left[
	\setlength\arraycolsep{6pt}
	\begin{matrix}
		q^{-2m}, &\quad bq^{-m}, &\quad cq^{-m} \\
		\multicolumn{3}{c}{
			\begin{matrix}
				q^{-m+1}/b, &q^{-m+1}/c 	
			\end{matrix}}
	\end{matrix} \hspace{2pt}
;q, \dfrac{q^{m+2}}{bc} \right] = \dfrac{(bq^{-m},cq^{-m};q)_m(q,bcq^{-2m};q)_{2m}}{(q,bcq^{-2m};q)_m(bq^{-m},cq^{-m};q)_{2m}}.
\end{equation} 
We now have  
\[ {}_3\phi_2 \left[
	\setlength\arraycolsep{6pt}
	\begin{matrix}
		q^{-2m}, &\quad bq^{-m}, &\quad cq^{-m} \\
		\multicolumn{3}{c}{
			\begin{matrix}
				q^{-m+1}/b, &q^{-m+1}/c 	
			\end{matrix}}
	\end{matrix} \hspace{2pt}
;q, \dfrac{q^{m+1}}{bc} \right] \]
\begin{align*}
&=\sum_{k=0}^{\infty}\dfrac{(q^{-2m},bq^{-m},cq^{-m};q)_k}{(q,q^{-m+1}/b,q^{-m+1}/c;q)_k}\left(\dfrac{q^{m+1}}{bc}\right)^k\\
&=\sum_{k=0}^{2m}\dfrac{(q^{-2m},bq^{-m},cq^{-m};q)_k}{(q,q^{-m+1}/b,q^{-m+1}/c;q)_k}\left(\dfrac{q^{m+1}}{bc}\right)^k\quad (\text{since}\,(q^{-2m})_k=0\,\text{for all}\,k>2m)\\
&=\sum_{k=0}^{2m}\dfrac{(q^{-2m},bq^{-m},cq^{-m};q)_{2m-k}}{(q,q^{-m+1}/b,q^{-m+1}/c;q)_{2m-k}}\left(\dfrac{q^{m+1}}{bc}\right)^{2m-k}\\
&=\dfrac{(q^{-2m},bq^{-m},cq^{-m};q)_{2m}(q^{m+1}/bc)^{2m}}{(q,q^{-m+1}/b,q^{-m+1}/c;q)_{2m}}\sum_{k=0}^{2m}\dfrac{(q^{-2m},bq^{-m},cq^{-m};q)_k}{(q,q^{-m+1}/b,q^{-m+1}/c;q)_k}\left(\dfrac{q^{m+2}}{bc}\right)^k\numberthis\label{eq42}\\
&=\dfrac{(q^{-2m},bq^{-m},cq^{-m};q)_{2m}(q^{m+1}/bc)^{2m}}{(q,q^{-m+1}/b,q^{-m+1}/c;q)_{2m}}\sum_{k=0}^{\infty}\dfrac{(q^{-2m},bq^{-m},cq^{-m};q)_k}{(q,q^{-m+1}/b,q^{-m+1}/c;q)_k}\left(\dfrac{q^{m+2}}{bc}\right)^k\\
&=\dfrac{(q^{-2m},bq^{-m},cq^{-m},q,bcq^{-2m};q)_{2m}(bq^{-m},cq^{-m};q)_m(q^{m+1}/bc)^{2m}}{(q,q^{-m+1}/b,q^{-m+1}/c,bq^{-m},cq^{-m};q)_{2m}(q,bcq^{-2m};q)_m}\numberthis\label{eq43}
\end{align*}
\\
where (\ref{eq42}) follows using (\ref{eq14}) with appropriate substitutions and (\ref{eq43}) follows from (\ref{eq41}).
\\\par Firstly, we note that the series on the left-hand side of (\ref{eq17}) is an analytic function of $1/d$ provided $\left\lvert q^2/bcd\right\rvert<\left\lvert q/bcd\right\rvert<1$. If we set $1/d=q^m$ for any positive integer $m$ in (\ref{eq17}), we get
\[{}_3\psi_3 \left[
   \begin{matrix}
      b, &c, &q^{-m} \\
      q/b, &q/c, &q^{m+1}
   \end{matrix}
;q, \dfrac{q^{m+1}}{bc} \right]\]
\begin{align*}
&=\sum_{k=-\infty}^{\infty}\dfrac{(b,c,q^{-m};q)_k}{(q/b,q/c,q^{m+1};q)_k}\left(\dfrac{q^{m+1}}{bc}\right)^k\\
&=\sum_{k=-m}^{\infty}\dfrac{(b,c,q^{-m};q)_k}{(q/b,q/c,q^{m+1};q)_k}\left(\dfrac{q^{m+1}}{bc}\right)^k\quad (\text{since}\,1/(q^{m+1})_k=0\,\text{for all}\,k<-m)\\
&=\sum_{k=0}^{\infty}\dfrac{(b,c,q^{-m};q)_{k-m}}{(q/b,q/c,q^{m+1};q)_{k-m}}\left(\dfrac{q^{m+1}}{bc}\right)^{k-m}\\
&=\dfrac{(b,c,q^{m};q)_{-m}(q^{m+1}/bc)^{-m}}{(q/b,q/c,q^{m+1};q)_{-m}}\sum_{k=0}^{\infty}\dfrac{(q^{-2m},bq^{-m},cq^{-m};q)_k}{(q,q^{-m+1}/b,q^{-m+1}/c;q)_k}\left(\dfrac{q^{m+1}}{bc}\right)^k\\
&=\dfrac{(b,c,q^{-m};q)_{-m}(q^{-2m},bq^{-m},cq^{-m},q,bcq^{-2m};q)_{2m}(bq^{-m},cq^{-m};q)_m(q^{m+1}/bc)^{m}}{(q/b,q/c,q^{m+1};q)_{-m}(q,q^{-m+1}/b,q^{-m+1}/c,bq^{-m},cq^{-m};q)_{2m}(q,bcq^{-2m};q)_m}
\end{align*}
where the last equality above follows from (\ref{eq43}). Then simplifying the last expression above using (\ref{eq11}), (\ref{eq12}) and (\ref{eq13}) with appropriate substitutions, we get
\\
\[{}_3\psi_3 \left[
   \begin{matrix}
      b, &c, &q^{-m} \\
      q/b, &q/c, &q^{m+1}
   \end{matrix}
;q, \dfrac{q^{m+1}}{bc} \right] = \dfrac{(q,q/bc,q^{m+1}/b,q^{m+1}/c;q)_{\infty}}{(q/b,q/c,q^{m+1},q^{m+1}/bc;q)_{\infty}}.\]
\\
\par Thus, the two sides of (\ref{eq17}) constitute analytic functions of $1/d$ provided $\left\lvert q^2/bcd\right\rvert<\left\lvert q/bcd\right\rvert<1$ where we note that the first of these inequalities always holds simply because $\lvert q\rvert<1$ and the second inequality can be rearranged to give $\lvert 1/d\rvert<\lvert bc/q\rvert$ which is a disk of radius $\lvert bc/q\rvert$ centred about $0$. Thus, both the sides of (\ref{eq17}) agree on an infinite sequence of points $(q^m)_{m\in \mathbb{N}}$ which converges to the limit $0$ inside the disk $\left\{1/d\in\mathbb{C}:\lvert 1/d\rvert<\lvert bc/q\rvert\right\}$. Hence, (\ref{eq17}) is valid in general. This completes the proof of Theorem \ref{3psi31}.
\end{proof}
\bigskip

\subsection{Proof of Theorem \ref{3psi32}}
\begin{proof}
 Replacing $n$ by $2m+1$, $z$ by $q^2$, $a$ by $bq^{-m-1}$ and $b$ by $cq^{-m-1}$ in (\ref{eq113}), we get
\begin{equation}\label{eq44}
{}_3\phi_2 \left[
	\setlength\arraycolsep{6pt}
	\begin{matrix}
		q^{-2m-1}, &bq^{-m-1}, &cq^{-m-1} \\
		\multicolumn{3}{c}{
			\begin{matrix}
				q^{-m+1}/b, &q^{-m+1}/c 	
			\end{matrix}}
	\end{matrix} \hspace{2pt}
;q, \dfrac{q^{m+4}}{bc} \right] \\ =\dfrac{(-1)^m(q^{-2m-1})_{2m}(q^2/bc)_{m}q^{\frac{m^2+5m}{2}}}{(q^2)_{m-1}(q^{-m+1}/b,q^{-m+1}/c;q)_{m}}.
\end{equation}
We now have  
\[ {}_3\phi_2 \left[
	\setlength\arraycolsep{6pt}
	\begin{matrix}
		q^{-2m-1}, &bq^{-m-1}, &cq^{-m-1} \\
		\multicolumn{3}{c}{
			\begin{matrix}
				q^{-m+1}/b, &q^{-m+1}/c 	
			\end{matrix}}
	\end{matrix} \hspace{2pt}
;q, \dfrac{q^{m+2}}{bc} \right] \]
\begin{align*}
&=\sum_{k=0}^{\infty}\dfrac{(q^{-2m-1},bq^{-m-1},cq^{-m-1};q)_k}{(q,q^{-m+1}/b,q^{-m+1}/c;q)_k}\left(\dfrac{q^{m+2}}{bc}\right)^k\\
&=\sum_{k=0}^{2m+1}\dfrac{(q^{-2m-1},bq^{-m-1},cq^{-m-1};q)_k}{(q,q^{-m+1}/b,q^{-m+1}/c;q)_k}\left(\dfrac{q^{m+2}}{bc}\right)^k\quad (\text{since}\,(q^{-2m-1})_k=0\,\text{for all}\,k>2m+1)\\
&=\sum_{k=0}^{2m+1}\dfrac{(q^{-2m-1},bq^{-m-1},cq^{-m-1};q)_{2m+1-k}}{(q,q^{-m+1}/b,q^{-m+1}/c;q)_{2m+1-k}}\left(\dfrac{q^{m+2}}{bc}\right)^{2m+1-k}\\
&=\dfrac{(q^{-2m-1},bq^{-m-1},cq^{-m-1};q)_{2m+1}q^{2m^2+5m+2}}{(q,q^{-m+1}/b,q^{-m+1}/c;q)_{2m+1}(bc)^{2m+1}}\sum_{k=0}^{2m+1}\dfrac{(q^{-2m-1},bq^{-m-1},cq^{-m-1};q)_k}{(q,q^{-m+1}/b,q^{-m+1}/c;q)_k}\left(\dfrac{q^{m+4}}{bc}\right)^k\numberthis\label{eq45}\\
&=\dfrac{(q^{-2m-1},bq^{-m-1},cq^{-m-1};q)_{2m+1}q^{2m^2+5m+2}}{(q,q^{-m+1}/b,q^{-m+1}/c;q)_{2m+1}(bc)^{2m+1}}\sum_{k=0}^{\infty}\dfrac{(q^{-2m-1},bq^{-m-1},cq^{-m-1};q)_k}{(q,q^{-m+1}/b,q^{-m+1}/c;q)_k}\left(\dfrac{q^{m+4}}{bc}\right)^k\\
&=\dfrac{(-1)^m(q^{-2m-1},bq^{-m-1},cq^{-m-1};q)_{2m+1}(q^{-2m-1})_{2m}(q^2/bc)_mq^{\frac{5m^2+15m+4}{2}}}{(q,q^{-m+1}/b,q^{-m+1}/c;q)_{2m+1}(q^2)_{m-1}(q^{-m+1}/b,q^{-m+1}/c;q)_m(bc)^{2m+1}}\numberthis\label{eq46}
\end{align*}
\\
where (\ref{eq45}) follows using (\ref{eq14}) with appropriate substitutions and (\ref{eq46}) follows from (\ref{eq44}).
\\\par Firstly, we note that series on the left-hand side of (\ref{eq18}) is an analytic function of $1/d$ provided $\left\lvert q^4/bcd\right\rvert<\left\lvert q^2/bcd\right\rvert<1$. If we set $1/d=q^m$ for any positive integer $m$ in (\ref{eq18}), we get
\[{}_3\psi_3 \left[
   \begin{matrix}
      b, &c, &q^{-m} \\
      q^2/b, &q^2/c, &q^{m+2}
   \end{matrix}
;q, \dfrac{q^{m+2}}{bc} \right]\]
\begin{align*}
&=\sum_{k=-\infty}^{\infty}\dfrac{(b,c,q^{-m};q)_k}{(q^2/b,q^2/c,q^{m+2};q)_k}\left(\dfrac{q^{m+2}}{bc}\right)^k\\
&=\sum_{k=-m-1}^{\infty}\dfrac{(b,c,q^{-m};q)_k}{(q^2/b,q^2/c,q^{m+2};q)_k}\left(\dfrac{q^{m+2}}{bc}\right)^k\quad (\text{since}\,1/(q^{m+2})_k=0\,\text{for all}\,k<-m-1)\\
&=\sum_{k=0}^{\infty}\dfrac{(b,c,q^{-m};q)_{k-m-1}}{(q^2/b,q^2/c,q^{m+2};q)_{k-m-1}}\left(\dfrac{q^{m+2}}{bc}\right)^{k-m-1}\\
&=\dfrac{(b,c,q^{m};q)_{-m-1}(q^{m+2}/bc)^{-m-1}}{(q^2/b,q^2/c,q^{m+2};q)_{-m-1}}\sum_{k=0}^{\infty}\dfrac{(q^{-2m-1},bq^{-m-1},cq^{-m-1};q)_k}{(q,q^{-m+1}/b,q^{-m+1}/c;q)_k}\left(\dfrac{q^{m+2}}{bc}\right)^k\\
&=\dfrac{(-bc)^{-m}(b,c,q^{-m};q)_{-m-1}(q^{-2m-1},bq^{-m-1},cq^{-m-1};q)_{2m+1}(q^{-2m-1})_{2m}(q^2/bc)_mq^{\frac{3m^2+9m}{2}}}{(q^2/b,q^2/c,q^{m+2};q)_{-m-1}(q,q^{-m+1}/b,q^{-m+1}/c;q)_{2m+1}(q^2)_{m-1}(q^{-m+1}/b,q^{-m+1}/c;q)_m}
\end{align*}
\\
where the last equality above follows from (\ref{eq46}). Then simplifying the last expression above using (\ref{eq11}), (\ref{eq12}) and (\ref{eq13}) with appropriate substitutions, we get
\\
\[{}_3\psi_3 \left[
   \begin{matrix}
      b, &c, &q^{-m} \\
      q^2/b, &q^2/c, &q^{m+2}
   \end{matrix}
;q, \dfrac{q^{m+2}}{bc} \right] = \dfrac{(q,q^2/bc,q^{m+2}/b,q^{m+2}/c;q)_{\infty}}{(q^2/b,q^2/c,q^{m+2},q^{m+2}/bc;q)_{\infty}}.\]
\\
\par Thus, the two sides of (\ref{eq18}) constitute analytic functions of $1/d$ provided $\left\lvert q^4/bcd\right\rvert<\left\lvert q^2/bcd\right\rvert<1$ where we note that the first of these inequalities always holds simply because $\lvert q\rvert<1$ and the second inequality can be rearranged to give $\lvert 1/d\rvert<\lvert bc/q^2\rvert$ which is a disk of radius $\lvert bc/q^2\rvert$ centred about $0$. Thus, both the sides of (\ref{eq18}) agree on an infinite sequence of points $(q^m)_{m\in \mathbb{N}}$ which converges to the limit $0$ inside the disk $\left\{1/d\in\mathbb{C}:\lvert 1/d\rvert<\lvert bc/q^2\rvert\right\}$. Hence, (\ref{eq18}) is valid in general. This completes the proof of Theorem \ref{3psi32}.
\end{proof}

\section*{Acknowledgments}
The author would like to thank Alexander Berkovich for encouraging him to prove Theorems \ref{5psi51}, \ref{5psi52}, \ref{3psi31}, \ref{3psi32} and for his very helpful comments and suggestions. In addition, the author would like to thank George E. Andrews and Jonathan Bradley-Thrush for previewing a preliminary draft of this paper and for their helpful comments. The author would also like to thank the anonymous referee and the editor for their helpful comments.

\end{document}